%% file: main.tex
\documentclass[12pt,twoside,a4paper]{article}
\usepackage{listings}
\usepackage{xcolor}
\usepackage[utf8]{inputenc}
\usepackage{mathtools}
\usepackage[none]{hyphenat}
\usepackage{mathtools}

\input{preamble}

\definecolor{codegreen}{rgb}{0,0.6,0}
\definecolor{codegray}{rgb}{0.5,0.5,0.5}
\definecolor{codepurple}{rgb}{0.58,0,0.82}
\definecolor{backcolor}{rgb}{0.95,0.95,0.92}

\lstdefinestyle{mystyle}{
  backgroundcolor=\color{backcolor},   commentstyle=\color{orange},
  keywordstyle=\color{magenta},
  numberstyle=\tiny\color{codegray},
  stringstyle=\color{codepurple},
  basicstyle=\ttfamily\footnotesize,
  breakatwhitespace=false,         
  breaklines=true,                 
  captionpos=b,                    
  keepspaces=true,                 
  numbers=left,                    
  numbersep=5pt,                  
  showspaces=false,                
  showstringspaces=false,
  showtabs=false,                  
  tabsize=2
}

\title{\sc Equitable Dominator Coloring of Graphs}

\author{{\bf Phebe Sarah George$^1$, S.Madhumitha$^2$ and Sudev Naduvath$^3$}}
\affil{Department of Mathematics\\ CHRIST University, Bangalore, India.\\
$^1${\tt phebe.george@res.christuniversity.in} \\
$^2${\tt s.madhumitha@res.christuniversity.in} \\
$^3${\tt sudev.nk@christuniversity.in}}
\date{}

\begin{document}
\maketitle
\hrule

\begin{abstract}
This paper introduces a variant of domination-related coloring of graphs, called the equitable dominator coloring of graphs, which is a combination of equitable coloring and dominator coloring.
The minimum number of colors used in an equitable dominator coloring of a graph is its equitable dominator chromatic number. The equitable dominator coloring and the corresponding equitable dominator chromatic number  of some standard graph classes are investigated in this paper.

\par

 
\end{abstract}

\maketitle 
\bigskip
%
\section{Introduction}

For basic terminology in graph theory, we refer to \cite{harary1,dbw1} and for topics in graph coloring, refer to \cite{cz1,mk1}. Unless mentioned otherwise, all graphs discussed in this paper are simple, undirected, finite, and connected.

\noi \emph{Graph coloring} is an assignment of colors to the graph's vertices, edges, or faces.
A \emph{vertex coloring} of a graph $G$ is a mapping $c:V(G) \to\mathcal{C}$, where $\mathcal{C}=\{c_1,c_2,\dots, c_k\}$, is a set of colors. A \emph{proper vertex coloring} of $G$ is when no two adjacent vertices are assigned the same color and the minimum number of colors required in this coloring of $G$ is the \emph{chromatic number} of $G$, denoted by $\chi(G)$.
The set of all vertices assigned the color $c_i$ in a coloring $c$ is called a \emph{color class}, denoted by $V_i$. A set $v \in V(G)$ is said to \emph{dominate} a set  $S\subseteq V(G)$, if $v$ is adjacent to every element of $S$. 

A proper coloring of a graph $G$ in which the cardinalities of any two color classes differ by at most 1 is said to be an \emph{equitable coloring} of $G$ (see \cite{1973equitable}) and the minimum number of colors used in this coloring is called the \emph{equitable chromatic number} of $G$, denoted by $\chi_{e}(G)$. An extensive study on the topic of equitable coloring can be found in literature in \cite{knesar,mycielskian,hypergraphs,corona,products, trees,  bipartiteeq}. 

The \emph{dominator coloring} of a graph $G$ is a proper coloring of $G$ such that every vertex in $V(G)$ dominates at least one color class, possibly its own color class (see \cite{DRsafeclique}). The minimum number of color classes in this coloring is called the \emph{dominator chromatic number} of $G$, denoted by $\chi_{d}(G)$.
For $S\subseteq V(G)$, a vertex $v\in V(G)$ is called a \emph{dominator} of $S$ if $v$ dominates $S$. The  dominator coloring of trees, bipartite graphs, Petersen graph and various other graph classes were studied in \cite{treesdom,bipartite,petersen,product,classes}. 
Motivated by the above studies, a variant of domination-related coloring, called equitable dominator coloring of graphs is introduced and studied in this paper. 

\section{Equitable Dominator Coloring of Graphs}

\noindent The notion of equitable dominator coloring of a graph is defined as follows.

\begin{definition}\label{eqdom}{\rm
An \textit{equitable dominator coloring} of a graph $G$ is a proper coloring of $G$ such that every vertex in $V(G)$ dominates at least one color class, possibly its own color class and the cardinalities of the color classes differ by at most one.
The minimum number of colors used in an equitable dominator coloring of $G$ is the \emph{equitable dominator chromatic number} of $G$, and we denote it  by $\chi_{ed}(G)$.}
\end{definition}

\begin{figure}[h]
\centering
\begin{subfigure}{0.45\textwidth}
   \begin{tikzpicture}[scale=0.55] 
\vertex (1) at (180:1) {$c_1$};
\vertex (2) at (0:1) []{$c_2$};
\vertex (3) at (0:3) []{$c_1$};
\vertex (4) at (0:5) []{$c_2$};
\vertex (5) at (120:2.5) []{$c_2$};
\vertex (6) at (150:2.5) []{$c_2$};
\vertex (7) at (180:2.5) []{$c_2$};
\vertex (8) at (210:2.5) []{$c_2$};
\vertex (9) at (240:2.5) {$c_2$};
\vertex (10) at (10:6.5) []{$c_{1}$};
\vertex (11) at (350:6.5) []{$c_{1}$};
\path
(5) edge (1)
(6) edge (1)
(7) edge (1)
(8) edge (1)
(9) edge (1)
(2) edge (1)
(2) edge (3)
(3) edge (4)
(4) edge (10)
(4) edge (11);
\end{tikzpicture}
    \subcaption{Minimum proper coloring of $G$.}
\end{subfigure}%
\hfill
\begin{subfigure}{0.45\textwidth}
    \begin{tikzpicture}[scale=0.55]
    \vertex (1) at (180:1) {$c_1$};
\vertex (2) at (0:1) []{$c_2$};
\vertex (3) at (0:3) []{$c_3$};
\vertex (4) at (0:5) []{$c_2$};
\vertex (5) at (120:2.5) []{$c_2$};
\vertex (6) at (150:2.5) []{$c_2$};
\vertex (7) at (180:2.5) []{$c_3$};
\vertex (8) at (210:2.5) []{$c_3$};
\vertex (9) at (240:2.5) {$c_3$};
\vertex (10) at (10:6.5) []{$c_{1}$};
\vertex (11) at (350:6.5) []{$c_{1}$};
\path
(5) edge (1)
(6) edge (1)
(7) edge (1)
(8) edge (1)
(9) edge (1)
(2) edge (1)
(2) edge (3)
(3) edge (4)
(4) edge (10)
(4) edge (11);
    \end{tikzpicture}
    \subcaption{Minimum equitable coloring of $G$.}
\end{subfigure}
\par\vspace{1cm} 
\begin{subfigure}{0.45\textwidth}
    \begin{tikzpicture}[scale=0.55]
   \vertex (1) at (180:1) {$c_1$};
\vertex (2) at (0:1) []{$c_3$};
\vertex (3) at (0:3) []{$c_4$};
\vertex (4) at (0:5) []{$c_2$};
\vertex (5) at (120:2.5) []{$c_3$};
\vertex (6) at (150:2.5) []{$c_3$};
\vertex (7) at (180:2.5) []{$c_3$};
\vertex (8) at (210:2.5) []{$c_3$};
\vertex (9) at (240:2.5) {$c_3$};
\vertex (10) at (10:6.5) []{$c_{4}$};
\vertex (11) at (350:6.5) []{$c_{4}$};
\path
(5) edge (1)
(6) edge (1)
(7) edge (1)
(8) edge (1)
(9) edge (1)
(2) edge (1)
(2) edge (3)
(3) edge (4)
(4) edge (10)
(4) edge (11);
    \end{tikzpicture}
    \subcaption{Minimum dominator coloring of $G$.}
\end{subfigure}%
\hfill
\begin{subfigure}{0.45\textwidth}
    \begin{tikzpicture}[scale=0.55]
    \vertex (1) at (180:1) {$c_1$};
\vertex (2) at (0:1) []{$c_2$};
\vertex (3) at (0:3) []{$c_3$};
\vertex (4) at (0:5) []{$c_4$};
\vertex (5) at (120:2.5) []{$c_2$};
\vertex (6) at (150:2.5) []{$c_5$};
\vertex (7) at (180:2.5) []{$c_5$};
\vertex (8) at (210:2.5) []{$c_6$};
\vertex (9) at (240:2.5) {$c_6$};
\vertex (10) at (10:6.5) []{$c_{3}$};
\vertex (11) at (350:6.5) []{$c_{7}$};
\path
(5) edge (1)
(6) edge (1)
(7) edge (1)
(8) edge (1)
(9) edge (1)
(2) edge (1)
(2) edge (3)
(3) edge (4)
(4) edge (10)
(4) edge (11);
    \end{tikzpicture}
    \subcaption{Minimum equitable dominator coloring of $G$.}
\end{subfigure}
\caption{A graph $G$ with $\chi(G)< \chi_e(G)< \chi_d(G) < \chi_{ed}(G)$.}
\label{fig:figures}
\end{figure}

Based on the definitions of proper coloring, equitable coloring, dominator coloring and equitable dominator coloring of graphs, it follows that

\begin{enumerate}[label=(\roman*),left=1cm,itemsep=2mm]
    \item $\chi(G) \le \chi_{e}(G) \le \chi_{ed}(G)$,
    \item $\chi_{d}(G) \leq \chi_{ed}(G)$.
\end{enumerate}

The Figure \ref{fig:figures} gives an example of a graph $G$ where the inequalities are sharp. The following propositions gives the conditions when the inequalities are strict. Note that a vertex $v$ in a graph $G$ of order $n$ is said to be a \emph{universal vertex}  of $G$, if $deg(v)=n-1$.

\begin{proposition}\label{universal}
If $\Delta(G)=n-1$, then $\chi_e(G)=\chi_{ed}(G)$.
\end{proposition}

\begin{proof}
In any proper coloring of a graph $G$ with at least one universal vertex $v$, every vertex of $G$ dominates the color class $\{v\}$. Also, in every equitable coloring of $G$, as the cardinalities of the color classes differ by at most 1, the result follows.
\end{proof}

The converse of Proposition \ref{universal} does not hold, as we identify a family of graphs without a universal vertex for which $\chi_{ed}(G)=\chi_e(G)$, in the following theorem.

\begin{proposition}\label{completemulti}
$\chi_{ed}(K_{a_1,a_2,\ldots,a_s})=\chi(K_{a_1,a_2,\ldots,a_s})$, if $|{a_i-a_j}|\leq 1; 1 \le i \neq j \le s$.
\end{proposition}

\begin{proof}
Any minimum proper coloring of $K_{a_1,a_2,\ldots,a_s}$ with $|{a_i-a_j}|\leq 1; 1 \le i \neq j \le s$, is its minimum equitable coloring. Furthermore, every vertex of the graph dominates $s-1$ color classes. Thus the result follows.
\end{proof}

The converse of Proposition \ref{completemulti} does not hold as we can see that $\chi_{ed}(C_5)=\chi(C_5)=3$. For a  graph $G$ of order $n$, $\chi_{ed}(G)=n$ if and only if $G$ is either a $K_n$ or  $\overline{K}_n$. In the following proposition, we characterise the graphs for which $\chi_{ed}(G)=2$.

\begin{proposition}
For a graph $G$, $\chi_{ed}(G)=2$ if and only if $G=K_{a,b}$ such that $|a-b| \leq 1$.
\end{proposition}

\begin{proof}
In the case when $G=K_{a,b}$, it follows that $\chi_{ed}(G)=2$ from Proposition \ref{completemulti}. To prove the converse, let $\chi_{ed}(G)=2$ for some graph $G$. Since $\chi_{ed}(G)=2$,  there is an independent set of vertices assigned the color $c_1$, say $V_1$, and an independent set of vertices assigned color $c_2$, say $V_2$. Also, every vertex of $V_1$ is adjacent to every vertex of $V_2$ and vice-versa. In order to satisfy the condition of equitability $b=a-1, a,a+1$ and this concludes the result.
\end{proof}

\begin{theorem}\label{realisation}
For any $j \in \N$, there exists at least one graph $G$ such that $\chi_{ed}(G)-\chi_{d}(G)=j$. 
\end{theorem}

\begin{proof}
Consider the graph $G=K_{2,2}$. We know that the dominator chromatic number for any complete bipartite graph is 2 and by Proposition 2.3 $\chi_{ed}(G)=2$. 

Now consider the graph $K_{2,4}$. Here, $|V_1|=2$ and hence the partite set $V_2$ can be partitioned into equitable parts with respect to $|V_1|$. 
As $\chi_{ed}(K_{2,4})=3$, on adding three vertices to $V_2$ in each iteration and making them adjacent to all the vertices of $V_1$, we get a complete bipartite graph $K_{2,4+3i}$, with equitable dominator chromatic number $3+i; 1 \le i \le n-3$. This proves the result.
\end{proof}

\begin{figure}[h]
\centering
\begin{tikzpicture}[scale=0.55] 
\vertex (1) at (3.75,5) []{$c_1$};
\vertex (2) at (8.25,5) []{$c_1$};
\vertex (3) at (1,1) []{$c_2$};
\vertex (4) at (3,1) []{$c_2$};
\vertex (5) at (5,1) []{$c_3$};
\vertex (6) at (7,1) []{$c_3$};
\vertex (7)[dashed] at (9,1) []{$c_4$};
\vertex (8)[dashed] at (11,1) []{$c_4$};
\vertex (9)[dashed] at (13,1) []{$c_4$};
\path
(1) edge (3)
(1) edge (4)
(1) edge (5)
(1) edge (6)
(1) edge [dashed] (7)
(1) edge [dashed] (8)
(1) edge [dashed] (9)
(2) edge (3)
(2) edge (4)
(2) edge (5)
(2) edge (6)
(2) edge [dashed] (7)
(2) edge [dashed] (8)
(2) edge [dashed] (9);
\end{tikzpicture}
\caption{Construction of graph $G$ such that $\chi_{ed}(G)-\chi_{d}(G)=2$.}
\end{figure}

\begin{figure}[h]
\centering
\begin{tikzpicture}[scale=0.55] 
\vertex (1) at (3.75,5) []{$c_1$};
\vertex (2) at (8.25,5) []{$c_1$};
\vertex (3) at (1,1) []{$c_2$};
\vertex (4) at (3,1) []{$c_2$};
\vertex (5) at (5,1) []{$c_3$};
\vertex (6) at (7,1) []{$c_3$};
\vertex (7)[dashed] at (9,1) []{$c_4$};
\vertex (8)[dashed] at (11,1) []{$c_4$};
\vertex (9)[dashed] at (13,1) []{$c_4$};
\vertex (10)[dashed] at (15,1) []{$c_5$};
\vertex (11)[dashed] at (17,1) []{$c_5$};
\vertex (12)[dashed] at (19,1) []{$c_5$};
\path
(1) edge (3)
(1) edge (4)
(1) edge (5)
(1) edge (6)
(1) edge [dashed] (7)
(1) edge [dashed] (8)
(1) edge [dashed] (9)
(2) edge (3)
(2) edge (4)
(2) edge (5)
(2) edge (6)
(2) edge [dashed] (7)
(2) edge [dashed] (8)
(2) edge [dashed] (9)
(1) edge [dashed] (10)
(1) edge [dashed] (11)
(1) edge [dashed] (12)
(2) edge [dashed] (10)
(2) edge [dashed] (11)
(2) edge [dashed] (12);
\end{tikzpicture}
\caption{Construction of graph $G$ such that $\chi_{ed}(G)-\chi_{d}(G)=7$.}
\label{Figure3}
\end{figure}

A graph realisation mentioned in Theorem \ref{realisation} is illustrated in Figure \ref{Figure3}, in which the dotted vertices and edges represent the added vertices and edges based on given value of $j$.

\section{Equitable Dominator Chromatic Number for Certain Graph Classes}\label{section3}

\begin{theorem}\label{path}
For $n \geq 4$,
\begin{gather*}
\chi_{ed}(P_n)=
\begin{cases}
3\floor{\frac{n}{5}}+1, &  n \equiv 0,1 \pmod{5};\\
3\ceil{\frac{n}{5}}-1, &  n \equiv 2,3 \pmod{5};\\
3\ceil{\frac{n}{5}}, &  n \equiv 4 \pmod{5}.
\end{cases}
\end{gather*}
\end{theorem}

\begin{proof}
For $P_n:=v_1-v_2- \ldots-v_n$, consider the following coloring patterns.

\begin{case}
When $n \equiv 1 \pmod{5}$, let $c:V(P_n)\to\{c_1,c_2,\ldots\}$ be a coloring such that,
\begin{gather*}
c(v_{j})=
\begin{cases}
c_{3\floor{\frac{j}{5}}+k_1}, &  j=5k+k_1, k_1=1,2,3;\\
c(v_{j-2}), & j\equiv 0,4 \pmod{5}.
\end{cases}
\end{gather*}

With respect to this coloring $c$ of $P_n$, the vertices $v_j; j \equiv 2 \pmod{5}$, dominate the color class $\{v_{j-1}\}$ and the vertices $v_j; j \equiv 0 \pmod{5}$, dominate the color class $\{v_{j+1}\}$. Also, the vertices $v_j; j \equiv 1 \pmod{5}$, dominate their own color classes and the vertices $v_j; j \equiv 3,4 \pmod{5}$, dominate the color classes $\{v_{j-1},v_{j+1}\}$. As the cardinality of the color classes of the colors used in  $c$ is at most $2$, it is an equitable dominator coloring of $P_n$ with $3(\frac{n-1}{5})+1$ colors.

Assume that there exists a coloring $c^{*}$ of $P_n; n=5k+1$, using $3(\frac{n-1}{5})$ colors. Therefore, with respect to $c^{*}$, there are $k$ color classes of cardinality $1$ and $2(\frac{n-1}{5})$ color classes of cardinality $2$, because a pendant vertex of $P_n$ can either dominate its own color class or the color class of its adjacent vertex. Hence, there exists no color class having cardinality greater than $2$ in any equitable dominator coloring of $P_n$. Two vertices in the same color class are at least at a distance $2$, because if the vertices are at distance greater than $2$, then at least two vertices adjacent to them need to have unique colors such that they dominate their own color classes; yielding a contradiction. 
 
\end{case}
\begin{case}
Consider a coloring $c^{\prime}$ of $P_n; n\not\equiv 1 \pmod{5}$, such that $c'(v_1)=c'(v_n)=c_1$,  $c'(v_2)=c_2, c'(v_{n-1})=c_3$, and for  $3 \le j \le n-2$,    
\begin{gather*}
c'(v_{j})=
\begin{cases}
c_{3\floor{\frac{j}{5}}+k_1+1}, &  j=5k+k_1, k_1=2,3,4;\\
c'(v_{j-2}), & j\equiv 0,1 \pmod{5}.
\end{cases}
\end{gather*}
With respect to the coloring $c'$ of $P_n$, the vertices $v_1,v_2$ dominate the color class $\{v_2\}$ and the vertices $v_{n-1},v_n$ dominate the color class $\{v_{n-1}\}$. Also, the vertices $v_j; j \equiv 3 \pmod{5}$, dominate the color class $\{v_{j-1}\}$, the vertices $v_j; j \equiv 2 \pmod{5}$, dominate their own color classes, the vertices $v_j; j \equiv 0,4 \pmod{5}$ dominate the color class $\{v_{j-1},v_{j+1}\}$ and the vertices $v_j; j \equiv 1 \pmod{5}$, dominate the color class $\{v_{j+1}\}$.  As the cardinality of the color classes of the colors used in  $c'$ is at most $2$, it is an equitable dominator coloring of $P_n$, in this case. Also, owing to the arguments mentioned in Case 1, it can be established that $c'$ is a minimum equitable dominator coloring of $P_n$. 
\end{case}
\end{proof}

\begin{theorem}\label{cycle}
For $n \geq 3$,
\begin{gather*}
\chi_{ed}(C_n)=
\begin{cases}
3\floor{\frac{n}{5}}+3, &  n \equiv 4 \pmod{5};\\
3\floor{\frac{n}{5}}+r, &  n \equiv r \pmod{5}, 0 \le r \le 3.
\end{cases}
\end{gather*}
\end{theorem}

\begin{proof}
Consider a cycle $C_n:= v_1-v_2-\ldots-v_n-v_1$ whose vertices are assigned colors as follows. Note that $v_{n+j}=v_j$.

\textit{Case 1:-} Let $n \equiv 0,1,2 \pmod{5}$.
\begin{gather*}
c(v_{j})=
\begin{cases}
c_{3\floor{\frac{j}{5}}+r}, &  j=5k_1+r, r=0,1,2;\\
c(v_{j-2}), & j\equiv 3,4 \pmod{5}.\\
\end{cases}
\end{gather*}

Based on the coloring protocol given, the vertices $v_j; j \equiv 1 \pmod{3}$, dominate the color class of the color assigned to the vertex $v_{j-1}$. The vertices $v_j; j \equiv 0 \pmod{5}$, dominate their own color classes. The vertices $v_j; j \equiv 2,3 \pmod{5}$, dominate the color class of the color assigned to the vertex $v_{j-1}$ and the vertices $v_j; j \equiv 4 \pmod{5}$, dominate the color class $\{v_{j+1}\}$. The cardinality of every color class in this coloring is at most $2$; thus satisfying the condition of equitability. Hence the coloring is an equitable dominator coloring and the result follows in this case.

\textit{Case 2:-} Let $n \equiv 3,4 \pmod{5}$
The vertices $v_i; 1 \le i \le n-1$ are assigned colors as mentioned in Case 1 and follows the property of equitable dominator coloring as justified in Case 1. 
However, in this case, the vertex $v_n$ needs to be assigned an unique color since the vertex $v_1$ cannot dominate the color class assigned to the vertex $v_2$ or its own color class and hence the result follows.

Since the vertices $\{v_i: 1 \le i \le n-1\}$ forms a path $P_{n-1}$ the minimality of the coloring follows from the Proof of Theorem 3.1. The optimality condition follows in $C_n$ as $v_1 \sim v_n$, $v_n$ needs to be assigned a unique color so that the vertex $v_1$ satisfies the condition of dominator coloring, this concludes the result.
\end{proof}

A \emph{bi-star} $S_{a,b}$ is a graph obtained by joining the central vertices of two-star graphs $K_{1,a}$ and $K_{1,b}$ by an edge.

\begin{theorem}
For $2 \le a \le b $, $\chi_{ed}(S_{a,b})=2+\ceil{\frac{a+b}{2}}$.
\end{theorem}

\begin{proof}
Let $u,v$ to be the support vertices of $S_{a,b}$ and $u_i;1\leq i\leq a$, and $v_j;1\leq j\leq b$, to be the pendant vertices adjacent to $u$ and $v$, respectively. Define a coloring $c:V(S_{a,b})\to\{c_1,c_2,\ldots,c_{2+\ceil{\frac{a+b}{2}}}\}$ as follows. For a vertex $w\in V(S_{a,b})$, 
\begin{gather*}
c(w)=
\begin{cases}
c_1, & w=u;\\
c_2, &  w=v;\\
c_{i+2}, & w \in \{u_i,v_i; 1 \le i \le a\};\\
c_{a+2+\ceil{\frac{i}{2}}}, & w =v_{a+i}; 1 \le i \le b-a. 
\end{cases}
\end{gather*}

We observe that, with respect to $c$, the vertices in $\{u\} \cup \{u_i:1 \le i \le a\}$, dominate the color class $\{u\}$ and the vertices in $\{v\} \cup \{v_i:1 \le i \le b\}$ dominate the color class $\{v\}$. Here, the cardinality of all the color classes is at most $2$ and hence, $\chi_{ed}(S_{a,b}) \le 2+ \ceil{\frac{a+b}{2}}$. 

Assume there exists an equitable dominator coloring of $S_{a,b}$, say $c'$, with fewer colors. In this case, either three pendant vertices are assigned the same color or one of the support vertex is assigned the same color as that of a pendant vertex, leading to a contradiction as a pendant vertex can dominate only its own color class or the color class of its adjacent vertex. 
\end{proof}

\begin{figure}[h]
\centering
\begin{subfigure}[b]{0.5\linewidth}
\centering
\begin{tikzpicture}[scale=0.55] 
\vertex (1) at (180:1) {$c_1$};
\vertex (2) at (0:1) []{$c_2$};
\vertex (3) at (135:2.5) []{$c_6$};
\vertex (4) at (165:2.5) []{$c_5$};
\vertex (5) at (195:2.5) []{$c_4$};
\vertex (6) at (225:2.5) []{$c_3$};
\vertex (7) at (345:2.5) {$c_5$};
\vertex (8) at (315:2.5) []{$c_4$};
\vertex (9) at (285:2.5) []{$c_3$};
\vertex (10) at (15:2.5) []{$c_6$};
\vertex (11) at (45:2.5) []{$c_7$};
\vertex (12) at (75:2.5) []{$c_7$};
\path

(1) edge (2)
(1) edge (3)
(1) edge (4)
(1) edge (5)
(1) edge (6)
(2) edge (7)
(2) edge (8)
(2) edge (9)
(2) edge (10)
(2) edge (11)
(2) edge (12)

;
\end{tikzpicture}
\end{subfigure}
\quad
\begin{subfigure}[b]{0.3\linewidth}
\centering
\begin{tikzpicture}[scale=0.55] 
\vertex (1) at (180:1) {$c_1$};
\vertex (2) at (0:1) []{$c_2$};
\vertex (3) at (220:3.5) []{$c_3$};
\vertex (4) at (200:3.5) []{$c_4$};
\vertex (5) at (180:3.5) []{$c_5$};
\vertex (6) at (160:3.5) []{$c_6$};
\vertex (13) at (140:3.5) []{$c_7$};
\vertex (7) at (280:3.5) {$c_3$};
\vertex (8) at (300:3.5) []{$c_4$};
\vertex (9) at (320:3.5) []{$c_5$};
\vertex (10) at (340:3.5) []{$c_6$};
\vertex (11) at (0:3.5) []{$c_7$};
\vertex (12) at (20:3.5) []{$c_8$};
\vertex (14) at (40:3.5) []{$c_8$};
\vertex (15) at (60:3.5) []{$c_9$};
\vertex (16) at (80:3.5) []{$c_9$};
\path

(1) edge (2)
(1) edge (3)
(1) edge (4)
(1) edge (5)
(1) edge (6)
(1) edge (13)
(2) edge (7)
(2) edge (8)
(2) edge (9)
(2) edge (10)
(2) edge (11)
(2) edge (12)
(2) edge (14)
(2) edge (15)
(2) edge (16)
;
\end{tikzpicture}
\end{subfigure}
\quad
\begin{subfigure}[b]{0.45\linewidth}
\centering
\begin{tikzpicture}[scale=0.55] 
\vertex (1) at (180:1) {$c_1$};
\vertex (2) at (0:1) []{$c_2$};
\vertex (3) at (218.55:3.5) []{$c_3$};
\vertex (4) at (192.85:3.5) []{$c_4$};
\vertex (5) at (167.15:3.5) []{$c_5$};
\vertex (6) at (141.45:3.5) []{$c_6$};
\vertex (7) at (282.9:3.5) {$c_3$};
\vertex (8) at (308.6:3.5) []{$c_4$};
\vertex (9) at (334.3:3.5) []{$c_5$};
\vertex (10) at (0:3.5) []{$c_6$};
\vertex (11) at (25.7:3.5) []{$c_7$};
\vertex (12) at (51.4:3.5) []{$c_7$};
\vertex (14) at (77.1:3.5) []{$c_8$};
\path
(1) edge (2)
(1) edge (3)
(1) edge (4)
(1) edge (5)
(1) edge (6)
(2) edge (7)
(2) edge (8)
(2) edge (9)
(2) edge (10)
(2) edge (11)
(2) edge (12)
(2) edge (14)
;
\end{tikzpicture}
\end{subfigure}\caption{Equitable dominator coloring of some bi-stars.}
\end{figure}

A \emph{helm graph} denoted by $H_{1,t,t}$ is a graph obtained by joining a pendant edge to each vertex of degree $3$ of a wheel graph $W_{1,t}=K_1+C_t$. 
\section{Equitable Dominator Coloring of Complete Bipartite Graphs}

Note that any proper coloring of $K_{a,b}$ with two colors is its dominator coloring; whereas, the property of equitability is not satisfied unless $|a-b|\leq 1$. Therefore, we use the concept of equitable partition of an integer defined as follows to find the equitable dominator coloring of $K_{a.b}$.

 An equitable partition of an integer $n$ is such that the integer $n$ is expressed as the sum of one or more positive integers such that the integers differ by at most 1 (see \cite{equitableinteger}).

 Finding the equitable partition of an integer $b$ with respect to the equitable partition of an integer $a$ is complex. Hence, a Python program is developed in order to find the equitable dominator chromatic number of a complete bipartite graph $K_{a,b}$ \label{Bipartite}. 
 
\begin{lstlisting}[language=Python]
def integerpart(integer):
    partites = set()
    partites.add((integer, ))
    for x in range(1, integer):
        for y in integerpart(integer - x):
            partites.add(tuple(sorted((x, ) + y)))
    return partites

def equitable(partites):
    listed=list(partites)
    length=len(partites)
    L1=[]
    for i in range (length):
        element=listed[i]
        element1=list(element)
        result=all(elements-element1[0]==1 or 
        elements-element1[0]==-1 or elements-element1[0]==0
        for elements in element1)
        if (result):
            L1.append(element)
        else:
            None
    return (L1) 

a=int(input("Enter first integer:"))
A=integerpart(a)

B=equitable(A)
B1= sorted(B, key=lambda x: len(x))
print(B1)

c=int(input("Enter the second integer:"))
C=integerpart(c)

D=equitable(C)
D1= sorted(D, key=lambda x: len(x))
print(D1)

def bipartite(B1,D1):
    listed1=list(B1)
    listed2=list(D1)
    length1=len(B1)
    length2=len(D1)
    if length1>length2:
        length3=length2
    else:
        length3=length1
    
    Pairs=[]
    for i in range (length1):
        for j in range (length2):
            elm1=listed1[i]
            elm2=listed2[j]
            elmA=list(elm1)
            elmB=list(elm2)
            combined = elmA + elmB
            result = max(combined) - min(combined) <= 1
            if (result):
                Pairs.append((elm1,elm2))
            else:
                None
    return (Pairs)

L=bipartite(B1,D1)

num=L[0]

x=[j for i in num for j in i]
x1=len(x)
print(x1)

print('The equitable dominator chromatic integer for K_{a,b} 
where a is %s and j is %s is %s.'% (a,c,x1))
\end{lstlisting}

In the algorithm above, the function \verb|integerpart()| is defined to find all possible integer partitions of an integer $n$, which is saved as a set and returned.
Next, the  \verb|equitable()| returns only those partitions from the output of \verb|integerpart()| that are equitable by comparing each $i$-th element of the set with all the other $j$-th elements of the set and is appended and returned in a list $L1$, whose elements are tuples and the elements of each tuple are the equitable parts of the given integer $n$.
When the user inputs the value of $a$, all the partitions of $a$ are returned, further the equitable partitions of $a$ are generated. Similarly, the user inputs the value of $b$, from which only the equitable partitions of $b$ are obtained. 
The equitable partitions of $a,b$ is then sorted and saved in lists $B1, D1$. A function  \verb|bipartite()| is then defined in which a list  \verb|Pairs| generated by comparing each element of each tuple from the list generated for $a$ with each element of each tuple from the list generated for $b$. The \verb|Pairs| list which is sorted, and the tuple with the minimum number elements in the list \verb|L| is saved in the variable \verb|num| and finally, $x1$ gives us the most optimal way of partitioning $a$ and $b$ in an equitable manner.

\section{Equitable Dominator Chromatic Number of Graph Complements}

We begin by discussing the equitable dominator chromatic number of the complement of paths and cycles. Note that $\overline{P}_2=2K_1$, $\overline{P}_3=K_1 \cup K_2$, $\overline{P}_4=P_4$ and hence, $\chi_{ed}(\overline{P}_2)=2$, $\chi_{ed}(\overline{P}_3)=3$, $\chi_{ed}(\overline{P}_4)=2$. Also, we observe that $\overline{C}_3=3K_1$ and $\overline{C}_4$ is $2K_2$, for which, $\chi_{ed}(\overline{C}_3)=3$ and  $\chi_{ed}(\overline{C}_4)=4$. Therefore, we consider the following result for $n \geq 5$.

Furthermore, $\overline{P}_n$ and $\overline{C}_n$ is a $\ceil{\frac{n}{2}}$-partite graph with cardinality of each part at most $2$.

\begin{theorem}
For $n \geq 5$, $\chi_{ed}(\overline{P}_n)=\chi_{ed}(\overline{C}_n)=\ceil{\frac{n}{2}}$.
\end{theorem}

\begin{proof}

For $V(\overline{P}_n)=\{v_i: 1 \le i \le n\}$, consider a coloring $c:V(\overline{P}_n)\to \mathcal{C}$ such that $c(v_i)=c(v_{i+1})=c_{\ceil{\frac{i}{2}}}$, where $i \equiv 1 \pmod{2}$.  This coloring is an equitable dominator coloring, since $\overline{P}_n$ is a $\ceil{\frac{n}{2}}$-partite graph with cardinality of each part at most $2$. Hence, the equitable dominator chromatic number of $\overline{P}_n$ is $\ceil{\frac{n}{2}}$. Also, the above argument holds for cycles.
\end{proof}

Note that $\overline{W}_{1,t}= \overline{C}_t \cup K_1$ and hence the following result is immediate.

\begin{corollary}
 For $t\geq 5$, $\chi_{ed}(W_{1,t})=1+\chi_{ed}(\overline{C}_t)$.
\end{corollary}

\begin{theorem}
For $a,b \ge 2$, $\chi_{ed}(\overline{S}_{a,b})=a+b$. 
\end{theorem}

\begin{proof}
Let $u,v$ to be the support vertices of $S_{a,b}$ and $u_i;1\leq i\leq a$, and $v_j;1\leq j\leq b$, to be the pendant vertices adjacent to $u$ and $v$, respectively. The pendant vertices in $S_{a,b}$ forms a clique $K_{a+b}$ in $\overline{S}_{a,b}$ and hence $\chi_{ed}(\overline{S}_{a,b})\geq a+b$. In $\overline{S}_{a,b}$, $u$ (resp. $v$) being adjacent to all the $v_j;1\leq j\leq b$ (resp. $u_i;1\leq i\leq a$), $u$ (resp. $v$) can be assigned any of the colors assigned to any of the $u_i;1\leq i\leq a$ (resp. $v_j;1\leq j\leq b$), in any equitable coloring of $\overline{S}_{a,b}$.

This coloring satisfies the condition of dominator coloring since the vertices $\{u\} \cup \{u_i: 1 \le i \le a\}$ dominate the color class of the colors assigned to the vertices $\{v_j: 1 \le j \le b\}$, the vertices $\{v\} \cup \{v_i: 1 \le i \le b\}$, dominate the color class of the colors assigned to the vertices $\{u_i: 1 \le i \le a\}$. Hence, $\chi_{ed}(\overline{S}_{a,b})=a+b$.
\end{proof}

\section{Conclusions}

This paper introduces the notion of equitable dominator coloring and determines the corresponding parameter of equitable dominator chromatic number for some graph classes and their complements. The equitable dominator chromatic number for a complete bipartite graph is found using a Python program and characterizations of graphs having some specific equitable dominator chromatic number have also been done. 

Some further directions of studies on equitable dominator coloring are mentioned below.

\begin{itemize}
    \item To study the parameter for various graph operations like join of graphs, strong product, normal product, etc.
    \item To study the parameter for generalized Petersen graphs
    \item To extend the study of $\chi_{ed}(G)$ to power of graphs.
\end{itemize}

\subsection*{Acknowledgement}
The first author would like to acknowledge her gratitude to her fellow researcher Ms.Sabitha Jose for their valuable suggestions and guidance throughout the work. 
\bigskip

\nocite{*}
\bibliographystyle{abbrv}

\end{document}

%% file: preamble.tex
\usepackage{graphicx} 
\usepackage{etoolbox}

\usepackage{amssymb,amsmath,amsthm,amsfonts,latexsym} 
\usepackage{correctmathalign}
\usepackage[mathscr]{euscript}

\usepackage{mathtools}
\usepackage[none]{hyphenat}
\DeclarePairedDelimiter{\ceil}{\lceil}{\rceil}
\DeclarePairedDelimiter\floor{\lfloor}{\rfloor}

\makeatletter
\DeclareRobustCommand*{\bfseries}{%
  \not@math@alphabet\bfseries\mathbf
  \fontseries\bfdefault\selectfont
  \boldmath
}
\makeatother

\usepackage[labelsep=quad,indention=10pt,format=hang]{caption}
\usepackage{subcaption}

\usepackage[pdftex,bookmarks,colorlinks=false]{hyperref}

\usepackage{geometry}
\usepackage[auth-lg,affil-sl]{authblk}
\usepackage{float} 
\usepackage{dpfloat}
\usepackage{multirow} 
\usepackage{hhline} %
\usepackage{booktabs} 
\usepackage{longtable} 
\usepackage{diagbox} 

\usepackage{cite}
\usepackage{url}
\usepackage{breakcites} 
\usepackage{breqn} 
\allowdisplaybreaks

\usepackage{cite}

\usepackage{color}
\usepackage{colordvi}
\usepackage{colortab}
\usepackage{colortbl}
\usepackage{colorweb}

\usepackage{bullcntr}%
\usepackage{doi}

\usepackage{epic}
\usepackage{epigraph}

\usepackage{enumitem}

\usepackage{calc,tikz}
\usetikzlibrary{}
\usetikzlibrary{intersections, arrows, arrows.meta, automata, er, calc, mindmap, folding, positioning,patterns, decorations.markings, fit, 
shapes, matrix, positioning, shapes.geometric, through}
\tikzstyle{vertex}=[circle, draw, inner sep=1pt, minimum size=4pt]
\newcommand{\vertex}{\node[vertex]}

\tikzstyle{ann} = [fill=white,font=\footnotesize,inner sep=1pt]
\tikzstyle{arrow} = [thick,<-->,>=stealth]

\newcommand{\noi}{\noindent}
\newcommand{\N}{\mathbb{N}}

\newtheorem{theorem}{Theorem}[section]
\newtheorem{definition}[theorem]{Definition}

\newtheorem{proposition}[theorem]{Proposition}
\newtheorem{corollary}[theorem]{Corollary}

\newcommand{\keywordsname}{Keywords}

\newcommand{\mscname}{MSC 2020}

\definecolor{ududff}{rgb}{0.30196078431372547,0.30196078431372547,1}
\definecolor{xdxdff}{rgb}{0.49019607843137253,0.49019607843137253,1}

\newtheoremstyle{casesty}         
{}                   
{}                   
{\upshape}           
{}                   
{\itshape}         
{:-}                   
{1em}                
{}                   
\theoremstyle{casesty}
\newtheorem{case}{Case}

\newtheoremstyle{subcasesty}         
{}                   
{}                   
{\upshape}           
{}                   
{\itshape}          
{:-}                   
{1em}                
{}                   
\theoremstyle{subcasesty}

\newtheoremstyle{schemesty}         
{}                   
{}                   
{\upshape}           
{}                   
{\itshape}          
{:-}                   
{1em}                
{}                   
\theoremstyle{schemesty}